\documentstyle[art10]{article}
\title{Involutions of a canonical curve.}
\author{Luis Fuentes\thanks{Supported by an F.P.U.
fellowship of Spanish Government} 
\and Manuel Pedreira}
\date{}

\newtheorem{teo}{Theorem}[section]
\newtheorem{defin}[teo]{Definition}
\newtheorem{prop}[teo]{Proposition}

\newtheorem{rem}[teo]{Remark}

\newtheorem{teo2}{Theorem}[subsection]

\newtheorem{prop2}[teo2]{Proposition}
\newtheorem{cor2}[teo2]{Corollary}
\newtheorem{lemma2}[teo2]{Lemma}

\newtheorem{rem2}[teo]{Remark}

\def\mas{_1}

\def\menos{_2}

\def\g2{\pi}
\def\cosa{{}}

\def\ov{\overline}

\def\p{{\bf P}}

\def\Te{{\cal O}}
\def\R{{\cal B}}

\font\euf=eufm10 at 12pt

\def\b{\mbox{\euf b}}
\def\bb{\mbox{\euf b}}

\def\A{{\alpha}}

\def\P{{\bf P}}

\def\cuatro{_{\cal K}}
\def\K{{\cal K}}
\def\k{{\cal K}}

\def\qed{\hspace{\fill}$\rule{2mm}{2mm}$}

\def\lrw{{\longrightarrow}}

\def\sub{{\subset}}

\def\vhi{\varphi}

\begin{document}

\maketitle

{\footnotesize{\bf Authors' address:} Departamento de Algebra, Universidad de Santiago
de Compostela. $15706$ Santiago de Compostela. Galicia. Spain. e-mail: {\tt
pedreira@zmat.usc.es}; \\ {\tt luisfg@usc.es}\\
{\bf Abstract:} We give a geometrical characterization of the ideal of quadrics containing a
ca\-noni\-cal curve with an involution. This implies to study involutions of rational normal
scrolls and Veronese surfaces.
}

\vspace{0.1cm}

\bigskip

{\bf Introduction:} Let $C$ be a nonhyperelliptic smooth curve of genus $\pi$. An involution of
$C$ is an automorphism $\vhi:C\lrw C$ such that $\vhi^2=id$. It induces a double cover
$\gamma:C\lrw C/\vhi=X$, where $X$ is a smooth curve of genus $g$. We say that $C$ has an
involution of genus $g$. By Hurwitz formula, we know that $\pi\geq 2g-1$. It is well know that
a general smooth curve of genus $\pi\geq 3$ has not nontrivial automorphisms. In particular a
smooth curve with an involution is not generic.

In this paper we give a geometric characterization of the ideal of quadrics containing the
canonical model of a nonhyperelliptic curve $C$ with an involution. Let
$C_{\K}\subset \P^{\pi-1}$ be the canonical model of $C$. We see that an involution of
$C_{\K}$ is a harmonic involution; that is, it can be extended to $\P^{\pi-1}$.

An involution $\ov{\vhi}$ of $\P^n$ has two complementary spaces of base points $S\mas$ and
$S\menos$. Moreover, $\ov{\vhi}$ induces an involution $\ov{\vhi}^*$ in the space of quadrics
of  $\P^n$. This involution has two spaces of base points: the {\em base quadrics}, that is,
quadrics containing the spaces $S\mas$ and $S\menos$ and the {\em harmonic quadrics}, that is,
quadrics such that $S\mas$ and $S\menos$ are polar respect to them. A subspace $\Sigma\subset
\P(H^0(\Te_{P^n}(2)))$ is called a {\em base-harmonic system} respect to $S\mas,S\menos$ when
it is a fixed space of $\ov{\vhi}^*$. In this case $\Sigma_b$ and $\Sigma_h$ will denote the
base quadrics and the harmonic quadrics of $\Sigma$ respectively. 

We prove the following Theorem: \\ \\ {\bf Theorem \ref{fundamental}} {\em
\ \ Let
$C_{\K}\subset \P^{\pi-1}$ be the canonical curve of genus $\pi$, with $\pi>4$. If
$C_{\K}$ has an involution of genus $g$ then $\pi\geq 2g-1$ and the quadrics of
$\P(H^0(I_{C_{\K}}(2)))$ are a base-harmonic system respect to the base spaces $\P^{g-1}$,
$\P^{\pi-g-1}$ that contains $(g-1)(\pi-g-2)$ independent base quadrics. Conversely, these
conditions are sufficient to grant the existence of an involution, except when:
\begin{enumerate}

\item $\pi=6,g=2$ and $C_{\K}$ has a $g^2_5$; or

\item $\pi=2g$,$2g+1$ or $2g+2$ and $C_{\K}$ is trigonal.

\end{enumerate} }
 \qed

First we prove that the conditions are sufficient when the curve $C_{\K}$ is complete
intersection of quadrics. The Enriq\"ues-Babbage Theorem says that $C_{\K}$ is the complete
intersection of quadrics except when it is trigonal or $C$ is a quintic smooth curve. In these
cases the quadrics containing the canonical curve intersect in a rational normal scroll and
in the Veronese surface respectively.

In order to examine the special cases we make an study of the harmonic involutions of the
rational normal scrolls and the Veronese surface. We compute the number of base quadrics on
each case. From this calculus we obtain the Corollaries \ref{coro1} and \ref{coro2}: \\ \\
{\bf Corollary \ref{coro1}}{\em \ \ The unique involutions on a trigonal canonical curve of
genus
$\pi$, $\pi>4$ are of genus
$\frac{\pi}{2},\frac{\pi-1}{2}$ or $\frac{\pi-2}{2}$. \qed} \\ \\ {\bf Corollary
\ref{coro2}}{\em \ \ The unique involutions on a smooth quintic plane curve are of genus
$2$.\qed }  

Furthermore, we make a particular geometrical study of the canonical curves of genus $4$ with
an involution of genus $2$ and genus $1$. 

Note that to compute the number of independent base quadrics containing the canonical curves
with and involution we need the result about the projective normality  of the canonical scrolls
(see \cite{pedreira2},$\S5$). We will follow the notation of \cite{fuentes} and
\cite{hartshorne} to work with scrolls and ruled surfaces.

We thanks Lawrence Ein by his interest on this work during his visit to our Department on
November, 2001.

\bigskip

\section{Harmonic involutions.}\label{involutions}

\begin{defin}\label{definvolucion}

Let $X\subset \P^n$ be a projective variety. An isomorphism $\vhi:X\lrw X$ is called an
involution if $\vhi^2=Id$. Moreover, if $\vhi$ is the restriction of an involution
$\ov{\vhi}:\P^n\lrw \P^n$ then $\vhi$ is called a harmonic involution.

\end{defin}

\begin{prop}\label{caracterizaarmonicas}

Let $X\subset \P^n$ be a linearly normal projective variety. An involution $\vhi:X\lrw X$ is
harmonic if and only if $\vhi^*(X\cap H)\sim X\cap H$ for all hyperplane $H$.

\end{prop}
{\bf Proof:} If $\vhi$ is harmonic we clearly have that $\vhi^*(X\cap H)\sim X\cap H$.

Conversely, if $\vhi^*(X\cap H)\sim X\cap H$ for all hyperplane, then we have an involution
$\vhi^*:H^0(\Te_X(1))\lrw H^0(\Te_X(1))$ that makes the following diagram commutative:

$$
\setlength{\unitlength}{5mm}
\begin{picture}(15,4)

\put(2.6,3){\makebox(0,0){$\P(H^0(\Te_X(1))^{\vee})$}}
\put(3.6,0){\makebox(0,0){$X$}}
\put(12.5,3){\makebox(0,0){$\P(H^0(\Te_X(1))^{\vee})$}}
\put(10.5,0){\makebox(0,0){$X$}}
\put(3.6,0.7){\vector(0,1){1.5}}
\put(10.5,0.7){\vector(0,1){1.5}}
\put(5.5,3){\vector(1,0){4}}
\put(5.5,0){\vector(1,0){4}}
\put(7.5,3.5){\makebox(0,0){$\P({\vhi^*}^{\vee})$}}
\put(7.5,0.6){\makebox(0,0){$\vhi$}}

\end{picture}
$$
so $\vhi$ extends to $\P^n$. \qed

{\bf Examples:}

\begin{enumerate}

\item {\em Any involution of a rational curve $D_n\subset \P^n$ is harmonic.}

It is sufficient to note that $\vhi^*(D_n\cap H)$ has degree $n$; because $D_n$ is a
rational curve, it follows that $\vhi^*(D_n\cap H)\sim D_n\cap H$.

\item {\em Any involution of a canonical curve $C_{\K} \subset \P^{\pi-1}$ of genus $\pi$ is
harmonic.}

The linear system $|\vhi^*(C_{\K}\cap H)|$ has degree $2\pi-2$ and it has dimension $\pi$, so
it is the canonical linear system and $\vhi^*(C_{\K}\cap H)\sim C_{\K}\cap H$.

\item {\em Any involution of a normal rational scroll $R_{n-1}\subset \P^n$ with invariant
$e>0$ is harmonic.}

Let $S_e=\P(\Te_{P^1}\oplus \Te_{P^1}(-e))$ be the ruled surface associated to $R_{n-1}$. We
know that
$H\cap R_{n-1}\sim X_0+\b f$. Since
$\vhi$ is an isomorphism,
$\vhi^*(X_0)^2=X_0^2$ and
$\vhi^*(f)^2=f^2=0$. But, in $S_e$ with $e>0$, $X_0$ is the unique curve with negative
self-intersection and $f$ is the unique curve with self-intersection $0$. From
this $\vhi^*(X_0)\sim X_0$ and $\vhi^*(f)\sim f$. Then:
$$
\vhi^*(H\cap R_{n-1})\sim \vhi^*(X_0+\b f)\sim \vhi^*(X_0)+\b \vhi^*(f)\sim X_0+\b f\sim H\cap
R_{n-1}
$$

\end{enumerate}

\bigskip

We recall some basic facts about involutions in a projective space.

\smallskip

\begin{prop}\label{espaciosfijos}

Any involution $\ov{\vhi}$ of $\P^n$ has two complementary spaces $S\mas,S\menos \\ \subset
\P^n$ of base points. In this way, the image of a point $P$ is the fourth harmonic of $P$,
$l\cap S\mas$ and $l\cap S\menos$, where $l$ is the  unique line passing through $P$ verifying
$l\cap S\mas\neq
\emptyset$ and $l\cap S\menos\neq \emptyset$.

Conversely, any pair of complementary spaces $S\mas,S\menos$ of $\P^n$ defines an invol\-ution
of
$\P^n$.\qed

\end{prop}

\begin{rem}\label{expresionmatricial}
{
If we take a base of $\P^n$, $W=\{ P_1,\ldots,P_{k+1},P'_1,\ldots,P'_{k'+1} \}$ where
$\langle P_1,\ldots,P_{k+1} \rangle=S\mas$ and $\langle P'_1,\ldots,P'_{k'+1} \rangle=S\menos$ then
the involution $\vhi$ is given by the matrix:
$$
M_{\vhi}=\left(
\begin{array}{cc}
{Id}&{0}\\
{0}&{-Id}\\
\end{array}\right)
$$ }\qed

\end{rem}

\begin{defin}\label{espaciosbase}

Under the above assumptions, we say that $\vhi$ is harmonic respect to $S\mas$ and $S\menos$.
Moreover, $S\mas$ and $S\menos$ are called the base spaces of $\vhi$.

\end{defin}

\begin{rem}\label{cuadricas}

{\em The linear isomorphism $\ov{\vhi}$ induces an isomorphism:
$$
\ov{\vhi}^*:\P(H^0(\Te_{P^n}(2)))\lrw \P(H^0(\Te_{P^n}(2)))
$$
Because $\ov{\vhi}$ is an involution, $\ov{\vhi}^*$ is an involution too. Therefore, it has two
complementary spaces of base points. Taking coordinates respect to the base $W$, a quadric
$Q\subset \P^n$ has a matrix:
$$
M_{Q}=\left(
\begin{array}{cc}
{A}&{C}\\
{C^t}&{B}\\
\end{array}\right)
$$

We see that $\ov{\vhi}^*(Q)=Q$ if and only if $M_{\ov{\vhi}^*(Q)}=M_{\vhi}M_QM_{\vhi}=\lambda
M_Q$ for some $\lambda\neq 0$, if and only if $A=B=0$ or $C=0$. In the first case $Q$ is
called a harmonic quadric and in the second case $Q$ is called a base quadric. The set of
harmonic quadrics will be denoted by $\P(H^0(\Te_{P^n}(2))_h)$ and the set of base quadrics by
$\P(H^0(\Te_{P^n}(2))_b)$. We have that:

\begin{enumerate}

\item $Q$ is a harmonic quadric if and only if $P'^tM_QP=0$ for each $P\in S\mas,P'\in
S\menos$; that is, if $S\mas$ and $S\menos$ are polar respect to $Q$.

\item $Q$ is a base quadric if and only if $S\mas,S\menos\subset Q$.

\end{enumerate} } \qed

\end{rem}

\begin{defin}\label{sistemaarmonico}

Let $\Sigma\subset \P(H^0(\Te_{P^n}(2)))$ be a projective subspace $\Sigma=\P(V)$ and let
$\ov{\vhi}:\P^n\lrw\P^n$ a harmonic involution respect to two spaces $S\mas,S\menos$. $\Sigma$
is called a base-harmonic system respect to $S\mas,S\menos$ when it is a fixed space of
$\ov{\vhi}^*$. 

\end{defin}

\begin{rem}\label{notasistemaarmonico}

{\em If we denote $\Sigma_h=\Sigma\cap \P(H^0(\Te_{P^n}(2))_h)$ and $\Sigma_b=\Sigma\cap
\P(H^0(\Te_{P^n}( \\ 2))_b)$, we see that $\Sigma$ is a base-harmonic system if and only if
$\Sigma=\Sigma_h+\Sigma_b$; that is, $V$ has a base composed by harmonic and base
quadrics.}\qed

\end{rem}

\begin{prop}\label{equivalencia}

Let $X\subset \P^n$ be a projective variety. 

\begin{enumerate}

\item If $X$ has a harmonic involution $\vhi$ respect to two spaces $S\mas,S\menos$ then
the system $\P(H^0(I_X(2)))$ is base-harmonic system respect to $S\mas,S\menos$.

\item Suppose that  $X$ is the complete intersection of quadrics. If $\P(H^0(I_X(2)))$ is
a base-harmonic system respect to $S\mas,S\menos$, then $X$ has a harmonic involution respect
to two spaces $S\mas,S\menos$.

\end{enumerate}

\end{prop}
{\bf Proof:}

\begin{enumerate}

\item Let $Q\in \P(H^0(I_X(2)))$, that is, $X$ is contained on $Q$. It is sufficient to show
that
$\ov{\vhi}^*(Q)$ contains $X$. Since $\ov{\vhi}(X)=X$ and $X\subset Q$, $X=\ov{\vhi}(X)\subset
\ov{\vhi}(Q)$ and the conclusion follows.

\item Let $\ov{\vhi}:\P^n\lrw \P^n$ the harmonic involution defined by the spaces $S\mas,S\menos$. If
$X$ is the complete intersection of quadrics, then $X=Q_1\cap\ldots\cap Q_k$ where
$\{ Q_1,\ldots ,Q_k\}$ is a base of $\P(H^0(I_X(2)))$. If $\P(H^0(I_X(2)))$ is a base-harmonic
system respect to $S\mas,S\menos$ we can choose a base of fixed quadrics, so 
$$
\ov{\vhi}(X)=\ov{\vhi}(Q_1\cap\ldots\cap Q_k)=\ov{\vhi}(Q_1)\cap\ldots\cap\ov{\vhi}(Q_k)=
Q_1\cap\ldots\cap Q_k=X
$$
and we can restrict $\ov{\vhi}$ to $X$. \qed

\end{enumerate}

\begin{prop}\label{cuadricasbase}

Let $X\subset \P^n$ be a projective variety and let $\vhi:X\lrw X$ a harmonic involution of
$X$ respect to spaces $S\mas,S\menos$. 

Let $F=\overline{\{P\in \P^n/P\in \langle x,\vhi(x) \rangle,
x\in X\}}$ be the variety of lines joining points of $X$ related by the involution. Then:
$$
\begin{array}{rl}
{h^0(I_{X,P^n}(2))_b}&{=h^0(I_{F,P^n}(2))_b=}\\
{}&{=h^0(I_{F,P^n}(2))-h^0(I_{F\cap
S\mas,S\mas}(2))-h^0(I_{F\cap S\menos,S\menos}(2))}\\
\end{array}
$$

\end{prop}
{\bf Proof:} Let us first prove that $H^0(I_{X,P^n}(2))_b=H^0(I_{F,P^n}(2))_b$.

Since $X\subset F$, a quadric containing $F$ contains $X$ too, so
$H^0(I_{F,P^n}(2))_b\subset H^0(I_{X,P^n}(2))_b$.

Conversely, if $Q\in H^0(I_{X,P^n}(2))_b$, $X,S\mas,S\menos\subset Q$. Therefore, each line $l$ of
$F$ meets $Q$ in four points: $(x,\vhi(x),l\cap S\mas,l\cap S\menos)$ and then it is contained on
$Q$. Thus,
$F\subset Q$ and $H^0(I_{X,P^n}(2))_b\subset H^0(I_{F,P^n}(2))_b$.

Now, let us consider the following exact sequence:
$$
0\lrw H^0(I_{F\cup S\mas,P^n}(2))\lrw H^0(I_{F,P^n}(2))\stackrel{\A}{\lrw} H^0(I_{F\cap
S\mas,S\mas}(2))
$$
Let us see that $\A$ is a surjective map. Let $Q\mas\subset S\mas$ be a quadric which contains
$F\cap S\mas$. Taking the cone of $Q\mas$ over $S\menos$, we obtain a quadric of $\P^n$ that contains
the lines joining $Q\mas$ and $S\menos$ so it contains $F$.
From this, we deduce that:
$$
h^0(I_{F\cup S\mas,P^n}(2))=h^0(I_{F,P^n}(2))-h^0(I_{F\cap S\mas,S\mas}(2))
$$
Similarly, we have
$$
0\lrw H^0(I_{F\cup S\mas\cup S\menos,P^n}(2))\lrw H^0(I_{F\cup S\mas,P^n}(2))\stackrel{\beta}{\lrw}
H^0(I_{F\cap S\menos,S\menos}(2))
$$
where $H^0(I_{F\cup S\mas\cup S\menos,P^n}(2))=H^0(I_{F,P^n}(2))_b$ and $\beta$ is a surjective map.
Therefore:
$$
\begin{array}{rl}
{h^0(I_{F,P^n}(2))_b}&{=h^0(I_{F\cup S\mas,P^n}(2))-h^0(I_{F\cap S\menos,S\menos}(2))=}\\
{}&{=h^0(I_{F,P^n}(2))-h^0(I_{F\cap S\mas,S\mas}(2))-h^0(I_{F\cap S\menos,S\menos}(2))}\\
\end{array}
$$\qed

{\em We finish this section by computing the dimension of the spaces of harmonic and base
quadrics in an involution over a normal rational curve.} 

Let $D_n\subset \P^n$ be a rational normal curve of degree $n$ and let $\vhi:D_n\lrw D_n$ an
involution of degree $2$. The lines joining points related by the involution generate a
rational normal ruled surface $R_{n-1}\subset \P^n$. This ruled surface has two directrix
curves $C\mas,C\menos$ on the base spaces $S\mas,S\menos$ of the involution. 

We now (\cite{hartshorne}, IV, $2.17$ and $2.19$) that $R_{n-1}\cong \P(\Te_{P^1}\oplus
\Te_{P^1}(-e))$ for some
$e\geq 0$. In this way the divisor of hyperplane sections of $R_{n-1}$ is $H\sim X_0+(n-1+e)/2
f$. Moreover the curve
$D_n$ corresponds to a
$2$-secant curve on the surface, so $D_n\sim 2X_0+kf$. Since $deg(D_n)=D_n.H=n$, we obtain
$k=e+1$. Because $D_n$ is irreducible, $D_n.X_0\geq 0$ and we obtain $e\leq 1$. We
deduce that
$e=0$ if $n$ is odd and $e=1$ if $n$ is even.

In this way we see that if $n$ is odd, the directrix curves $C\mas,C\menos$ have degree $(n-1)/2$ and
if $n$ is even, the directrix curves have degree $(n-2)/2$ and $n/2$.

By applying Proposition \ref{cuadricasbase} we have that:
$$
h^0(I_{D_n}(2))_b=h^0(I_{R_{n-1}}(2))-h^0(I_{C\mas}(2))-h^0(I_{C\menos}(2))
$$
The dimension of theses spaces are well known. Thus, we obtain:
$$
h^0(I_{D_n}(2))_b=\left\{
\begin{array}{l}
{(\frac{n-1}{2})^2\mbox{ if $n$ is odd}}\\
{}\\
{\frac{n(n-2)}{4}\mbox{ if $n$ is even}}\\
\end{array}\right.
$$
 
Note that when $n$ is odd, $D_n\sim 2X_0+f$, so $D_n.X_0=1$ with $C\mas,C\menos\sim X_0$. Thus, the
base points of the involution of $D_n$ are $D_n\cap C\mas$ and $D_n\cap C\menos$. When $n$ is
even, $D_n\sim 2X_0+2f$, so $D_n.X_0=0$ and $D_n.(X_0+f)=2$ where $C\mas\sim X_0$ and $C\menos\sim
X_0+f$. In this case the base points of the involutions of $D_n$ are the two points of
$D_n\cap C\menos$. \qed

\bigskip

\section{Involutions of the canonical curve.}\label{involutionscanonica}

Let $C_{\K}$ be the canonical curve of genus $\pi$ and let $\vhi:C_{\K}\lrw C_{\K}$ be an
involution. We saw that it is a harmonic involution. We will use the results obtained in
\cite{pedreira2}. The scroll generated by the
involution is a canonical scroll
$R_{\bb}$. We call genus of the involution $\vhi$ to the genus of the ruled surface $R_{\bb}$.
Thus we have a $2:1$-morphism $\gamma:C_{\K}\lrw X$. $R_{\bb}$ has a canonical directrix curve
$\ov{X_0}$ of genus $g$, and a nonspecial curve
$\ov{X_1}$ with degree
$\pi-1$. They lie on disjoint spaces $\P^{g-1}$ and $P^{\pi-g-1}$. The involution of $C_{\K}$
has $2(\pi-1-2(g-1))$ base points that are the ramifications of $\gamma$. We denote them by
$\R$ and we know that $\R\sim C_{\K}\cap \ov{X_1}$.

Let us compute the number of base quadrics containing $C_{\K}$, that is, the dimension of
$H^0(I_{C_{\K}}(2))_b$. By Proposition \ref{cuadricasbase} we know that:

$$
h^0(I_{C_{\K}}(2))_b=h^0(I_{R_{\bb}}(2))-h^0(I_{\ov{X_0},P^{g-1}}(2))-h^0(I_{\ov{X_1},P^{\pi-g-1}}(2))
$$

We can compute this dimension because we know that (see \cite{pedreira2}):
$$
\begin{array}{l}
{h^0(I_{R_{\bb}}(2))=h^0(\Te_{P^{\g2-1}}(2))-h^0(\Te_{S_{\bb}}(2H))+dim(s(H,
H))}\\
{h^0(I_{\ov{X_0}}(2))=h^0(\Te_{P^{g-1}}(2))-h^0(\Te_X(2\k))+dim(s(\k,\k))}\\

{h^0(I_{\ov{X_1}}(2))=h^0(\Te_{P^{\g2-g-1}}(2))-h^0(\Te_X(2\b))+dim(s(\b,\b))}\\
{}\\
{dim(s(H,H))=dim(s(\k,\k))+dim(s(\b,\b))}\\
\end{array}
$$
From this, we obtain:
$$
h^0(I_{C_{\K}}(2))_b=h^0(\Te_{P^{\g2-1}}(2))_b-h^0(\Te_X(\b+\k))
$$
Thus the number of base
quadrics containing $C_{\K}$ is:
$$
h^0(C_{\K}(2))_b=(g-1)(\g2-g-2)
$$ \qed

\begin{teo}\label{involucioncanonica}

Let $C_{\K}\subset \P^{\pi-1}$ be a nonhyperelliptic canonical curve of genus $\pi$ that it is
complete intersection of quadrics. $C_{\K}$ has an involution of degree $2$ if and only if
$\P(H^0(I_{C_{\K}}(2)))$ is harmonic respect to two disjoint complementary subspaces of
dimensions
$k$ and $\pi-k-2$, with
$k\leq \pi-k-1$. Moreover the involution has genus $g=k+1$ if and only if
$h^0(I_{C_{\K}}(2))_b=(g-1)(\pi-g-2)$.

\end{teo}
{\bf Proof:} The first assertion follows from Proposition \ref{equivalencia}.

By the above discussion we know the number of base quadrics
$h^0(I_{C_{\K}}(2))_b=(g'-1)(\pi-g'-2)$ where $g'$ is the genus of the involution. This genus
is $k+1$ or $\pi-k-1$. Suppose that $(g'-1)(\pi-g'-2)=(g-1)(\pi-g-2)$ with $g=k+1$. If
$g'=\pi-k-1=\pi-g$ then we obtain $\pi=2g$ so the genus $g'$ of the involution is $k+1$. \qed

\begin{rem}\label{icompleta}

{\em It is well known that a nonhyperelliptic canonical curve of genus $\pi\geq 5$ is the
complete intersection of quadrics, except when it is trigonal or when $\pi=6$ and it has a
$g^2_5$. In the first case the intersection of the quadrics of $\P(H^0(I_{C_{\K}}(2)))$ is a
rational normal scroll; in the second case it is the Veronese surface of $\P^5$.} \qed

\end{rem}

\bigskip

\subsection{Involutions of the canonical curve of genus $\pi=4$.}\label{involuciones4}

Let $C\cuatro\subset \P^3$ be a canonical curve of genus $4$. It is well known that this curve is
the complete intersection of a quadric and a cubic surface (see \cite{hartshorne},IV,Example
$5.2.2$). Suppose that
$C\cuatro$ has an involution
$C\cuatro\lrw X$ where $X$ is a smooth curve of genus $g$. We have that $\pi-1\geq 2g-2$, so then
genus of $X$ can be $1$ or $2$ (if $g=0$ the curve $C\cuatro$ is hyperelliptic).

\begin{prop2}\label{involucion41}

A canonical curve of genus $4$ has an elliptic involution if and only if it is the complete
intersection of an elliptic cubic cone $S$ and a quadric that doesn't pass through the vertex
of
$S$. 

\end{prop2}
{\bf Proof:} Suppose that $C\cuatro$ has an elliptic involution. We know that the involution
generates a scroll $R$. In this case the directrix curve $X_0$ has degree $g-1=0$, so
it is a point. Then
$R$ is an elliptic cone. If $Q$ is the unique quadric that contains $C\cuatro$, then necessary
$C\cuatro=Q\cap R$. Moreover, we know that $C\cuatro\cap X_0=\emptyset$, so $Q$ does not pass
through the vertex of $R$. Conversely, if $C\cuatro=Q\cap R$ and $Q$ does not pass
through the vertex of $R$, the generators of $S$ provide an elliptic involution of $C\cuatro$.
\qed

\begin{teo2}\label{involucion42}

A canonical curve of genus $4$ has an involution of genus $2$ if and only if it is the
complete intersection of a quadric and a cubic surface which has a harmonic involution
respect to polar lines respect to the quadric.

\end{teo2}
{\bf Proof:}

\begin{enumerate}

\item Suppose that $C\cuatro=Q_2\cap Q_3$ where $Q_2$ is a quadric and $Q_3$ is a cubic
surface with a harmonic involution $\vhi:Q_3\lrw Q_3$. Let $l$ and $l'$ the base spaces.
Suppose that they are polar respect to $Q_2$. From this, $\vhi(Q_2)=Q_2$ so $C\cuatro$ has a harmonic involution. Because the base spaces have dimension $1$, the involution has genus $2$.

\item Suppose that $C\cuatro$ has a harmonic involution $\vhi$ of genus $2$. From this, the base
spaces are two lines $l$ and $l'$. Let $R\subset \P^3$ the ruled surface generated by the
involution. We know that $R$ contains $l$ and $l'$ with multiplicities $3$ and $2$
respectively. Moreover,
$l'\cap C\cuatro=\emptyset$ and $l\cap C\cuatro$ consist of two points. We know that
$\P(H^0(I_{C\cuatro}(2)))=\{ Q_2
\}$ and $dim\P(H^0(I_{C\cuatro}(3)))=4$. In this way $C\cuatro=Q_2\cap Q_3$ for any cubic surface $Q_3\in
\P(H^0(I_{C\cuatro}(3)))$. Consider the involution $\ov{\vhi}:\P(H^0(I_{C\cuatro}(3)))\lrw
\P(H^0(I_{C\cuatro}(3)))$ induced by
$\vhi$. Let us see that there is a fixed irreducible element. Let $V$ be the set of reducible
elements of $\P(H^0(I_{C\cuatro}(3)))$. Because $Q_2$ is the unique quadric containing $C\cuatro$ then
$V=\{Q_2+H; H\subset \P^3 \}$ and $dim(V)=3$. Since the base points of $\ov{\vhi}$ generate
$\P(H^0(I_{C\cuatro}(3)))$, then there exist at least a fixed irreducible element. \qed

\end{enumerate}

\bigskip

\subsection{Involutions of the Veronese surface.}\label{involucionesV}

Let $v_{2,n}$ be the Veronese map of $\P^n$:
$$
\begin{array}{rcl}
{v_{2,n}:{\P^n}^*}&{\lrw}&{V_{2,n}\subset \P(H^0(\Te_{P^n}(2)))}\\
{[x_0:\ldots :x_n]}&{\lrw}&{[x_0^2:x_0x_1:\ldots :x_n^2]}\\
\end{array}
$$
We will denote the image of this map by $V_{2,n}$. If $n=2$, it is the Veronese surface and we
will denote it by $V_2$.

\begin{prop2}\label{iveronesevariedad}

The involutions of the Veronese variety $V_{2,n}$ are harmonic respect to two base
spaces which are the harmonic and base quadrics of two subspaces $S\mas,S\menos$ of $\P^n$.

\end{prop2}
{\bf Proof:} A harmonic involution $\vhi:\P^n\lrw \P^n$ respect to spaces $S\mas,S\menos$
induces a harmonic involution $\ov{\vhi}$ in $\P(H^0(\Te_{P^n}(2)))$ respect to
$\P(H^0(\Te_{P^n}(2))_h)$ and $\P(H^0(\Te_{P^n}(2))_b)$.

Because $\ov{\vhi}\circ v_{2,n}=v_{2,n}\circ \ov{\vhi}$, we see that $\ov{\vhi}$ restricts to
$V_{2,n}$. In this way we have a harmonic involution $\ov{\vhi}:V_{2,n}\lrw V_{2,n}$ respect
to the spaces
$\P(H^0(\Te_{P^n}(2))_h)$ and $\P(H^0(\Te_{P^n}(2))_b)$.

Conversely, given an involution $\eta$ of $V_2$. Applying the isomorphism $v_{2,n}$, we obtain
an involution $\vhi^*$ of ${\P^n}^*$ with base spaces ${S\mas}^*,{S\menos}^*$. In this way, the
dual map $\vhi=\vhi^{**}$ provides an involution $\ov{\vhi}$ of $\P(H^0(\Te_{P^n}(2)))$ such
that
$\ov{\vhi}|_{V_{2,n}}=\eta$.\qed

We will denote the Veronese surface by $V_2$.

\begin{cor2}\label{iveronesesuperficie}

Any nontrivial involution of the Veronese surface $V_2$ is harmonic respect to a line which
corresponds to the conics of $\P^2$ passing through a point $P$ and a line $r$, and a
$3$-dimensional space $V$ corresponding to the polar conics of $\P^2$ respect to $P$ and $r$.

\end{cor2}
{\bf Proof:} It is sufficient to note that an involution $\eta$ of $V_2$ is induced by and
involution $\vhi$ of $P^2$. If $\vhi$ is nontrivial, then its base spaces are a point $P$ and
a line $r$.\qed

\begin{prop2}\label{veroneseqbase}

Let $\eta$ be a nontrivial harmonic involution of the Veronese surface. Then
$h^0(\Te_{V_2}(2))_b=2$.

\end{prop2}
{\bf Proof:} Consider the Veronese map of $\P^2$:
$$
\begin{array}{rcl}
{v_2:{\P^n}^*}&{\lrw}&{V_2\subset \P(H^0(\Te_{P^n}(2)))}\\
{[x_0:x_1:x_2]}&{\lrw}&{[x_0^2:x_0x_1:\ldots :x_2^2]=[y_0:y_1:\ldots:y_5]}\\
\end{array}
$$

Let $Y$ be the matrix:
$$
\left(\begin{array}{ccc}
{y_0}&{y_1}&{y_2}\\
{y_1}&{y_3}&{y_4}\\
{y_2}&{y_4}&{y_5}\\
\end{array}\right)
$$
We know that $H^0(\Te_{V_2}(2))$ is generated by the quadrics of $\P^5$ whose equations are
defined by the minors of order $2$ of the matrix $Y$.

Moreover the base spaces of $\vhi$ are a line $l$ corresponding to the conics of $\P^2$
passing through a line $r$ and a point $P$ and a space $V$ corresponding to the polar conics
respect to $P$ and $r$.

Taking an adequate system of coordinates we can consider $P$ generated by the equations
$\{x_1=x_2=0\}$ and $r$ generated by the equation ${x_0=0}$.

A conic containing $P$ and $r$ has an equation $ax_0x_1+bx_0x_2=0$. Then the equations of $l$
are $\{y_0=y_3=y_4=y_5=0\}$.

A polar conic respect to $P$ and $r$ has an equation $ax_0^2+bx_1^2+cx_1x_2+dx_2^2=0$. Then
the equations of $V$ are $\{y_1=y_2=0\}$.

Applying the conditions to contain $V$ and $l$ to the equations of $H^0(\Te_{V_2}(2))$, we
obtain that the quadrics containing $V$ and $l$ are generated by
$\{y_1y_4-y_2y_3=y_1y_5-y_2y_4=0\}$. From this, $h^0(\Te_{V_2}(2))_b=2$. \qed

\bigskip

\subsection{Harmonic involutions of the rational ruled surfaces.}\label{involucionesR}

Let $R_{n-1}\subset \P^n$ a rational normal ruled surface of degree $n-1$, with $n>3$. Let
$\vhi:R_{n-1}\lrw R_{n-1}$ a harmonic involution of the surface. Then $\vhi$ conserves the
degree of the curves. From this, it applies generators into generators if $n>3$.
In this way, we have the following induced harmonic involutions:
$$
\begin{array}{rccl}
{\vhi_0:\P^1}&{\lrw}&{\P^1}&{}\\
{P}&{\lrw}&{Q}&{/\vhi(Pf)=Qf}\\
\end{array}
$$
where $\P^1$ parameterizes the generators.
$$
\begin{array}{rcl}
{\vhi_l:|D_l|}&{\lrw}&{|D_l|}\\
{D}&{\lrw}&{\vhi(D)}\\
\end{array}
$$
where $|D_l|$ is the linear system of curves of degree $l$.

Let $k$ be the degree of the curve of minimum degree of $R_{n-1}$:
\begin{enumerate}

\item If $k=\frac{n-1}{2}$, then there is a $1$-dimensional family of irreducible curves of
degree $k$. The involution $\vhi_k$ has at least $2$ base points, so there are two disjoint
curves $D_k$ that are invariant by $\vhi$.

\item If $k<\frac{n-1}{2}$, then there is a unique curve of minimum degree, so it is invariant
by $\vhi$. Moreover, if $l=n-1-k$ the linear system $|D_l|$ has dimension $l-k$. Its generic
curve is an irreducible curve disjoint from $D_k$. In particular, the set of reducible curves
of $|D_l|$ are an hyperplane composed by curves of the form $D_k+\sum f_i$. Thus, we have a harmonic involution:
$$
\begin{array}{rcl}
{\vhi_l:\P^{l-k}\cong |D_l|}&{\lrw}&{|D_l|\cong \P^{l-k}}\\
\end{array}
$$
We know that $\vhi_l$ has two disjoint spaces of base points. Both of them can not be contained
on the hyperplane (because they generate $\P^{n-k}$), so necessary there exists an irreducible
curve in $|D_l|$ that is fixed by the involution $\vhi_l$; that is, it is invariant by $\vhi$.

\end{enumerate}

We conclude the following proposition:

\begin{prop2}\label{conclusion1}

Given a harmonic involution on a rational normal ruled surface of degree $n-1$, there exist
two disjoint rational normal curves $D_k$,$D_l$ with degrees $k$ and $l$,
$k+l=n-1$, that are invariant by the involution.

\end{prop2}

Let $D_k\subset \P^k$ and $D_l\subset \P^l$ be the two invariant curves. The involution $\vhi$
restricts to these spaces. Thus, we have a harmonic involution $\vhi_k:\P^k\lrw \P^k$. It has
two invariant spaces $\P^{k\mas}$, $\P^{k\menos}$ with $k\mas+k\menos+1=k$. Similarly, the harmonic
involution $\vhi_l:\P^l\lrw \P^l$ has two invariant spaces $P^{l\mas}$, $P^{l\menos}$ with
$l\mas+l\menos+1=l$. From this, we have two possibilities for the base spaces $S\mas$,$S\menos$
of the involution $\vhi$:
$$
\begin{array}{ll}
{\begin{array}{l}
	{S\mas=\langle \P^{k\mas},\P^{l\mas} \rangle=\P^{k\mas+l\mas+1}}\\
	{S\menos=\langle \P^{k\menos},\P^{l\menos} \rangle=\P^{k\menos+l\menos+1}}\\
\end{array}}&
{\begin{array}{l}
	{S\mas=\langle \P^{k\mas},\P^{l\menos} \rangle=\P^{k\mas+l\menos+1}}\\
	{S\menos=\langle \P^{k\menos},\P^{l\mas} \rangle=\P^{k\menos+l\mas+1}}\\
\end{array}}\\
\end{array}
$$

Conversely if we have two harmonic involutions in $\P^k$ and $\P^l$ we can recuperate an
involution in $\P^n$. Note that this involution is not unique, because we have two ways to
define it. Moreover, in order to restrict the involution to the ruled surface $q:R_{n-1}\lrw
\P^1$ we need that the involutions in $\P^k$ and $\P^l$ restrict to $D_k$ and $D_l$ and that
they are compatible, that is, the images of the points on the same generator lay on the same
generator: $q(\vhi_k(D_k\cap P f))=q(\vhi_l(D_l\cap P f)), \forall P\in \P^1$.

Thus, if $\vhi_k$ and $\vhi_l$ verify these conditions we have a harmonic involution $\vhi$
that restricts to $R_{n-1}$.

\begin{prop2}\label{conclusion2}

A harmonic involution on a normal rational ruled surface $R_{n-1}$ defines two harmonic
involutions $\vhi_k$, $\vhi_l$ on two disjoint rational curves $D_k$, $D_l$ that generate the
surface. Moreover, they make commutative the diagram $(1)$.

Conversely, if two harmonic involutions $\vhi_k$, $\vhi_l$ on two rational curves generating a
rational ruled surface $R_{n-1}$ verifying $q(\vhi_k(D_k\cap P f))=q(\vhi_l(D_l\cap P f)),
\forall P\in \P^1$. , then they define two
possible harmonic involutions on $R_{n-1}$, taking the space bases generate by the space bases
of $\vhi_k$ and $\vhi_l$.

\end{prop2}

\begin{rem2}\label{compatibilidad}

{\em In order to obtain compatible involutions $\vhi_k$, $\vhi_l$ it is sufficient to define a
involution $\eta$ on $\P^1$ and to translate it to $D_k$ and $D_l$: 
$$
\begin{array}{c}
{\vhi_k(D_k\cap P f):=D_k\cap \eta(P) f}\\
{\vhi_l(D_l\cap P f):=D_l\cap \eta(P) f}\\
\end{array}
$$ for all $P\in \P^1$.

Moreover, if $\vhi_k$ and $\vhi_l$ are compatible involutions and one of them is the identity,
then the other one is the identity too.}\qed
\end{rem2}

We saw how are the (nontrivial) harmonic involutions on a normal rational curve $D_m\sub \P^m$:
\begin{enumerate}

\item If $m=2\mu$ the involution is defined by two base spaces $\P^{\mu}$, $\P^{\mu-1}$ such
that
$\P^{\mu}\cap D_m=P\cup Q$ and $\P^{\mu-1}\cap D_m=\emptyset$ ($P,Q$ base points).

\item If $m=2\mu+1$ then involution is defined by two base spaces $\P^{\mu}_1$, $\P^{\mu}_2$
such that
$\P^{\mu}_1\cap D_m=P$ and $\P^{\mu}_2\cap D_m=Q$ ($P,Q$ base points).

\end{enumerate}

In both cases the involution generates a normal rational ruled surface $R_{m-1}$ of degree
$m-1$, whose directrix curves lay on the space bases. We call them base curves.

Thus, let $\vhi$ be a harmonic involution on $R_{n-1}$. Let $\vhi_k$, $\vhi_l$ be the
harmonic involutions induced on the directrix curves $D_k$, $D_l$. Let $\P^{k\mas}$, $\P^{k\menos}$,
$\P^{l\mas}$, $\P^{l\menos}$ be the base spaces of $\vhi_k$ and $\vhi_l$. We know that the base
spaces of $\vhi$ are $S\mas=\langle \P^{k\mas},\P^{l\mas} \rangle$, $S\menos=\langle
\P^{k\menos},\P^{l\menos} \rangle$. Let $C_{k\mas}$, $C_{k\menos}$, $C_{l\mas}$, $C_{l\menos}$ the corresponding
base curves. Let $F$ be the variety of lines that join the points of the involution:
$F=\overline{\{P\in \P^n/P\in \langle x,f(x)\rangle,x\in R_{n-1}\}}$. Let us identificate
$F\cap S\mas$ and $F\cap S\menos$.

\begin{lemma2}\label{variedadderectas}

The variety $F\cap S\mas$ $(F\cap S\menos)$ is a normal rational ruled surface of degree
$k\mas+l\mas=n\mas$ $(k\menos+l\menos=n\menos)$ generated by the directrix curves $C_{k\mas}$ and
$C_{l\mas}$ ($C_{k\menos}$ and
$C_{l\menos}$). We call it base ruled surface $R_{n\mas-1}$ ($R_{n\menos-1}$).

\end{lemma2}
{\bf Proof:} Given a point $P\in R_{n-1}$, consider the line $r=\langle P,f(P) \rangle$ of $F$.
$r$ meets
$S\mas$ in a point $P\mas$ that corresponds to project $P$ from $S\menos$ onto $S\mas$.
Thus, given a generator $f\in R_{n-1}$, the lines of $F$ defined by the points of $f$ meet
$S\mas$ in a line $f\mas$; this line is the projection of $f$ from $S_b$. Moreover, since
$f$ meets
$D_k$ and
$D_l$, its projection on $S\mas$ meets $C_{k\mas}$ and $C_{l\mas}$. In this way we see that the
generator of $R_{n-1}$ project into lines joining $C_{k\mas}$ and $C_{l\menos}$, so $F\cap S\mas$
is the rational ruled surface defined by these directrix curves. \qed

We saw that a harmonic involution on a normal rational ruled surface is defined by
the involutions of the directrix curves $D_k$ and $D_l$. From this, we distinguish several
types of involutions:

\begin{enumerate}

\renewcommand{\theenumii}{\arabic{enumii}}
\renewcommand{\theenumiii}{\arabic{enumiii}}
\renewcommand{\theenumiv}{\arabic{enumiv}}
\renewcommand{\labelenumii}{\theenumi.\theenumii.}
\renewcommand{\labelenumiii}{\theenumi.\theenumii.\theenumiii.}
\renewcommand{\labelenumiv}{\theenumi.\theenumii.\theenumiii.\theenumiv}

\makeatletter
\renewcommand{\p@enumiv}{\theenumi.\theenumii.\theenumiii.}
\makeatother

\item {\em $\vhi_k$ and $\vhi_l$ are the identity.}

Then the base spaces of $\vhi$ are the spaces $\P^k$ and $\P^l$ that contain the directrix
curves. All the generators are invariants by $\vhi$ and the variety $F$ is the ruled surface
$R_{n-1}$.

\item {\em $\vhi_k$ and $\vhi_l$ are not trivial.}

\begin{enumerate}

\item $n-1=2\lambda$ (even).

\begin{enumerate}

\item $k=2\mu, l=2(\lambda-\mu)$.

Then the involutions on $D_k$ and $D_l$ have the following base spaces and base curves:

$C_{\mu}\in \P^{\mu}$, $C_{\mu-1}\in \P^{\mu-1}$, with $P_k,Q_k\in C_{\mu}$ base points of
$D_k$.

$C_{\lambda-\mu}\in \P^{\lambda-\mu}$, $C_{\lambda-\mu-1}\in \P^{\lambda-\mu-1}$, with
$P_l,Q_l\in C_{\lambda-\mu}$ base points of $D_l$.

Then, the base spaces of $\vhi$ are:

\begin{enumerate}

\item Case A:

$S\mas=\langle \P^{\mu},\P^{\lambda-\mu}\rangle=\P^{\lambda+1}\ni P_k,Q_k,P_l,Q_l$.

$S\menos=\langle \P^{\mu-1},\P^{\lambda-\mu-1}\rangle=\P^{\lambda-1}$.

Where the generators $f_P,f_Q\in \P^{\lambda+1}$ are fixed.

\item Case B: \label{B1}

$S\mas=\langle \P^{\mu},\P^{\lambda-\mu-1}\rangle=\P^{\lambda}\ni P_k,Q_k$.

$S\menos=\langle \P^{\mu-1},\P^{\lambda-\mu}\rangle=\P^{\lambda}\ni P_l,Q_l$.

Where the generators $f_P,f_Q$ are invariant (not fixed).

\end{enumerate}

\item $k=2\mu+1, l=2(\lambda-\mu)-1$.

Then the involutions on $D_k$ and $D_l$ have the following base spaces and base curves:

$P_k\in C_{\mu}\in \P^{\mu}$, $Q_k\in C_{\mu}\in \P^{\mu}$, with $P_k,Q_k$ base
points of $D_k$.

$P_l\in C_{\lambda-\mu}\in \P^{\lambda-\mu}$, $Q_l\in C_{\lambda-\mu}\in \P^{\lambda-\mu}$,
with $P_l,Q_l$ base points of $D_l$.

Then, the base spaces of $\vhi$ are:

\begin{enumerate}

\item Case C:

$S\mas=\langle \P^{\mu},\P^{\lambda-\mu-1}\rangle=\P^{\lambda}\ni P_k,P_l$.

$S\menos=\langle \P^{\mu},\P^{\lambda-\mu-1}\rangle=\P^{\lambda}\ni Q_k,Q_l$.

Where the generators $f_P\in \P^a,f_Q\in S\menos$ are fixed.

\item Case B: (similar to case \ref{B1}).

$S\mas=\langle \P^{\mu},\P^{\lambda-\mu-1}\rangle=\P^{\lambda}\ni P_k,Q_l$.

$S\menos=\langle \P^{\mu},\P^{\lambda-\mu-1}\rangle=\P^{\lambda}\ni P_l,Q_k$.

Where the generators $f_P,f_Q$ are invariant (not fixed).

\end{enumerate}

\end{enumerate}

\item $n-1=2\lambda-1$ even.

Then the curves $D_k$ and $D_l$ have degrees $k=2\mu$ and $l=2(\lambda-\mu)-1$. The base
spaces and base curves are:

$C_{\mu}\in \P^{\mu}$, $C_{\mu-1}\in \P^{\mu-1}$, with $P_k,Q_k\in C_{\mu}$ base points of
$D_k$.

$P_l\in C_{\lambda-\mu-1}\in \P^{\lambda-\mu-1}$, $Q_l\in C_{\lambda-\mu-1}\in
\P^{\lambda-\mu-1}$, with $P_l,Q_l$ base points of $D_l$.

In any case, the base spaces of $\vhi$ are:

\begin{enumerate}

\item Case D:

$S\mas=\langle \P^{\mu},\P^{\lambda-\mu-1}\rangle=\P^{\lambda}\ni P_k,Q_k,P_l$.

$S\menos=\langle \P^{\mu},\P^{\lambda-\mu-1}\rangle=\P^{\lambda}\ni Q_l$.

Where $f_P$ is a fixed generator and $f_Q$ are an invariant generator.

\end{enumerate}

\end{enumerate}

\end{enumerate}

Let $\vhi:R_{n-1}\lrw R_{n-1}$ a harmonic involution over the rational normal ruled surface
$R_{n-1}\subset \P^n$. By the proposition \ref{equivalencia} we know that
$H^0(I_{R_{n-1}}(2))$ is a base-harmonic system; that is,
$H^0(I_{R_{n-1}}(2))=H^0(I_{R_{n-1}}(2))_h\oplus H^0(I_{R_{n-1}}(2))_b$. Let us see the
dimension of these spaces. We know that $h^0(I_{R_{n-1}}(2))=\left(^{n-1}_{\;\;\, 2}\right)$.
We will treat each case separated:

\begin{enumerate}

\item All generators are invariant by the involution.

We use the proposition \ref{cuadricasbase}. In this case $F=R_{n-1}$, $F\cap \P^a=D_k$ and
$F\cap \P^b=D_l$. Thus,
$$
h^0(I_{R_{n-1}}(2))_b=h^0(I_{R_{n-1}}(2))-h^0(I_{D_k}(2))-h^0(I_{D_l}(2))=kl
$$
and
$$
h^0(I_{R_{n-1}}(2))_h=h^0(I_{R_{n-1}}(2))-h^0(I_{R_{n-1}}(2))_b=\left(^{n-1}_{\;\;\,
2}\right)-kl
$$

\item The generic generator is not invariant by the involution.

We use the proposition \ref{cuadricasbase}. But in this case, $F\cap S\mas=R_{n\mas-1}$ and $F\cap
S\menos=R_{n\menos}$. Then we have:
$$
h^0(I_{R_{n-1}}(2))_b=h^0(I_{F}(2))-h^0(I_{R_{n\mas-1}}(2))-h^0(I_{R_{n\menos-1}}(2))
\eqno(2)
$$

A quadric containing $R_{n-1}\cup R_{n\mas-1}$ meets each line of $F$ in three points.
Then such quadric contains $F$, so $H^0(I_F(2))=H^0(I_{R_{n-1}\cup R_{n\mas-1}}(2))$. 
Consider the exact sequence:
$$
0\lrw H^0(I_{R_{n-1}\cup R_{n\mas-1}}(2))\lrw H^0(I_{R_{n-1}}(2))\stackrel{\alpha}{\lrw}
H^0(\Te_{R_{n\mas-1}}(2-Y))
$$
where $Y=R_{n-1}\cap R_{n\mas-1}$. Then:
$$
h^0(I_F(2))\geq h^0(I_{R_{n-1}}(2))-h^0(\Te_{R_{n\mas-1}}(2-R_{n-1}\cap R_{n\mas-1}))
$$
and applying $(2)$ we obtain in each case:

\begin{enumerate}

\renewcommand{\theenumii}{\Alph{enumii}}
\renewcommand{\labelenumii}{\theenumii.}

\item ($n-1=2\lambda, S\mas=\P^{\lambda-1},S\menos=\P^{\lambda+1},f_Q,f_P$ fixed generators.)
$$
h^0(I_{R_{n-1}}(2))_b\geq \lambda (\lambda-1)
$$

\item ($n-1=2\lambda, S\mas=\P^{\lambda}_1,S\menos=\P^{\lambda}_2,f_Q,f_P$ invariant (not fixed)
gen\-er\-ators, with $f_P\cap \P^{\lambda}_i=P_i$, $f_Q\cap P^{\lambda}_i=Q_i$.)
$$
h^0(I_{R_{n-1}}(2))_b\geq \lambda (\lambda-1)
$$

\item $(n-1=2\lambda, S\mas=\P^{\lambda}_1,S\menos=\P^{\lambda}_2,f_Q,f_P$ fixed
generators, with $f_P\in \P^{\lambda}_1$, $f_Q\in \P^{\lambda}_2$.)
$$
h^0(I_{R_{n-1}}(2))_b\geq \lambda (\lambda-1)+1
$$

\item $(n-1=2\lambda-1, S\mas=\P^{\lambda-1},S\menos=\P^{\lambda},f_Q$ fixed
generator in $\P^{\lambda}$, and $f_q$ invariant generator, with $f_Q\cap P^{\lambda-1}=Q_1$.)
$$
h^0(I_{R_{n-1}}(2))_b\geq (\lambda-1)^2
$$

\end{enumerate}

\end{enumerate}

Now, let us compute the harmonic quadrics.

Let $E_k$ the set of $k+1$ generic points in $D_k$ and $E_l$ the set of $l+1$ generic points
on $D_l$. Note that a harmonic quadric that passes through a point of $R_{n-1}$ passes
through the image point too. From this, a quadric passing through $E_k (E_l)$ meets $D_k (D_l)$
in $2k+2 (2l+2)$ points because $D_k (D_l)$ is invariant by the involution. Moreover, a
harmonic quadric that contains $D_k$ and $D_l$, contains the invariant (not fixed) generators
too, because they meet each space base in a point.

Finally, a harmonic quadric containing $D_k$ and $D_l$ and passing through $m$ generic points
of $R_{n-1}$ contains their images ($2m$ points) and the corresponding $2m$ generators. Let
$E_m$ be $m$ generic points of $R_{n-1}$ and let $E$ be $E_k\cup E_l \cup E_m$.

If $Q$ is a harmonic quadric passing through the points of $E$, then $D_k\cup
D_l\cup \{$invariant generators$\}\cup 2mf\subset Q\cap R_{n-1}$. If $2m>2(n-1)-(k+l)-$number
of invariant generators, then $R_{n-1}\subset Q$ and we have the exact sequence:
$$
0\lrw H^0(I_{R_{n-1}}(2))_h\lrw H^0(I_{\P^{n-1}}(2))_h\lrw H^0(\Te_E(2))
$$
From this:
$$
h^0(I_{R_{n-1}}(2))_h\geq h^0(\Te_{\P^{n-1}}(2))_h-(n+1+m)
$$

In each case we obtain:

\begin{enumerate}

\renewcommand{\theenumi}{\Alph{enumi}}

\item There are not invariant generators. Taking $m=\lambda+1$ we have:
$$
h^0(I_{R_{n-1}}(2))_h\geq \left(^{\lambda+1}_{\;\;\, 2}\right) + \left(^{\lambda+3}_{\;\;\,
2}\right)-(3\lambda+3)
$$

\item There are two invariant generators. Taking $m=\lambda$ we have:
$$
h^0(I_{R_{n-1}}(2))_h\geq \left(^{\lambda+2}_{\;\;\, 2}\right) + \left(^{\lambda+2}_{\;\;\,
2}\right)-(3\lambda+2)
$$

\item There are not invariant generators. Taking $m=\lambda+1$ we have:
$$
h^0(I_{R_{n-1}}(2))_h\geq \left(^{\lambda+2}_{\;\;\, 2}\right) + \left(^{\lambda+2}_{\;\;\,
2}\right)-(3\lambda+3)
$$

\item There is an invariant generator. Taking $m=\lambda$ we have:
$$
h^0(I_{R_{n-1}}(2))_h\geq \left(^{\lambda+1}_{\;\;\, 2}\right) + \left(^{\lambda+2}_{\;\;\,
2}\right)-(3\lambda+1)
$$

\end{enumerate}

We see that the sum of the bounds computed for the harmonic and base quadrics is the
quadrics of $H^0(I_{R_{n-1}}(2))$ in all cases, so these bounds are reached in all
cases, and we have the number of base quadrics:

\begin{prop2}\label{baseregladaracional}

Let $R_{n-1}\subset \P^n$ be a rational normal scroll of degree $n-1$. Let
$\vhi:R_{n-1}\lrw R_{n-1}$ be a harmonic involution. Then we have the following cases:

\begin{enumerate}

\renewcommand{\theenumii}{\arabic{enumii}}
\renewcommand{\theenumiii}{\Alph{enumiii}}

\item All the generators are invariant. There are two directrix curves of base point $D_k,D_l$
with $k+l=n-1$. They lay on the base spaces $\P^k$, $\P^l$.
$$
h^0(I_{R_{n-1}}(2))_b=kl
$$

\item There are two invariant (fixed or not):

\begin{enumerate}

\item $n-1=2\lambda$ ($n$ even).

\begin{enumerate}

\item The base spaces are $\P^{\lambda-1}$,$\P^{\lambda+1}$. There are two fixed generators in
$\P^{\lambda+1}$.

$$
h^0(I_{R_{n-1}}(2))_b=\lambda(\lambda-1)
$$

\item The base spaces are $\P^{\lambda}$,$\P^{\lambda}$. There is a fixed generator in each
of them.

$$
h^0(I_{R_{n-1}}(2))_b=\lambda(\lambda-1)+1
$$

\item The base spaces are $\P^{\lambda}$,$\P^{\lambda}$. There are not fixed generators.

$$
h^0(I_{R_{n-1}}(2))_b=\lambda(\lambda-1)
$$

\end{enumerate}

\item $n-1=2\lambda-1$ ($n$ odd).

\begin{enumerate}

\item[D.] The base spaces are $\P^{\lambda-1}$,$\P^{\lambda}$. There is a fixed generator in
$\P^{\lambda}$.

$$
h^0(I_{R_{n-1}}(2))_b=(\lambda-1)^2
$$

\end{enumerate}

\end{enumerate}

\end{enumerate}

\end{prop2} \qed

\bigskip

\subsection{Involutions of the canonical curve of genus $\pi>4$.}\label{involucionesC}

We have investigated all the possible cases where the quadrics that contain a canonical curve
are a base-harmonic system:

\begin{teo2}\label{uno}

The unique cases where the system of quadrics containing a canonical curve $C_{\K}$ of genus
$\pi$, $\pi>4$ are a base-harmonic system respect to base spaces $\P^{g-1}$, $\P^{\pi-g-1}$
($\pi\geq 2g-1>0
$) with $b$ independent base quadrics are:

\begin{enumerate}

\item $b=(g-1) (\pi-g-2)$

\begin{enumerate}

\item If $\pi=2g,2g+1,2g+2$ and the curve has a $g^1_3$ or an
involution of genus $g$, or both of them; except if $\pi=6$ and $g=2$ when $C_{\K}$ can
have a $g^2_5$.

\item If $\pi\neq 2g,2g+1,2g+2$ and the curve $C_{\K}$ has a $\gamma^1_2$ of genus
$g$.

\end{enumerate}

\item $b=(g-1)(\pi-g-2)+1$, $\pi=2g-1r$ or $\pi=2g$ and $C_{\K}$ has a $g^1_3$ (not a
$\gamma^1_2$).

\item $b=(g-1)(\pi-g-1)$ and $C_{\K}$ has a $g^1_3$ (not a $\gamma^1_2$).

\end{enumerate}

\end{teo2}

\begin{cor2}\label{coro1}

The unique involutions on a trigonal canonical curve of genus $\pi$, $\pi>4$ are of genus
$\frac{\pi}{2},\frac{\pi-1}{2}$ or $\frac{\pi-2}{2}$.

\end{cor2}
{\bf Proof:} Let $C_{\K}\subset \P^n$ be a trigonal canonical curve and let $R_{\pi-2}$ be the
ruled surface of trisecants. Suppose that $C_{\K}$ has an involution of genus $g$. Then the
system of quadrics
$\P(H^0(I_{C_{\K}}(2)))$ is a base-harmonic system with $(g-1)(\pi-g-2)$ independent base
quadrics. Since
$R_{\pi-2}=\bigcap_{Q\supset C_{\K}}Q$, we have a harmonic involution over $R_{\pi-2}$. By
Proposition \ref{baseregladaracional} we know that $\pi=2g$, $\pi=2g+1$ or $\pi=2g+2$ and the
conclusion follows. \qed

\begin{cor2}\label{coro2}

The unique involutions on a smooth quintic plane curve are of genus $2$. \qed

\end{cor2} 

\begin{teo2}\label{casoparticular}

Let $C_{\K}\subset \P^{\pi-1}$ a canonical curve of genus $\p$, $\pi>4$. Then $C_{\K}$ has
an involution of genus $1$ if and only if the quadrics of $\P(H^0(I_{C_{\K}}(2)))$ are a
base-harmonic system respect to a point and a space $\P^{\pi-2}$ without base quadrics.

\end{teo2}
{\bf Proof:} If $C_{\K}$ has an involution of genus $1$, we know that the involution generates
an elliptic cone, the system of quadrics $\P(H^0(I_{C_{\K}}(2)))$ is harmonic and it hasn't
base quadrics respect to $\P^0$ and $\P^{\pi-2}$.

Conversely, if the system of quadrics $\P(H^0(I_{C_{\K}}(2)))$ is harmonic respect to $\P^0$
and $\P^{\pi-2}$, necessary it hasn't base quadrics, because the quadric containing $C_{\K}$
are reducible. If $C_{\K}$ is not trigonal we have an involution of genus $1$ in $C_{\K}$.
If $C_{K}$ is trigonal, by Corollary \ref{coro1}, $\pi=2,3,4$. But we
have supposed that $\pi>4$. \qed

\begin{teo2}\label{fundamental}

Let $C_{\K}\subset \P^{\pi-1}$ be the canonical curve of genus $\pi$, with $\pi>4$. If
$C_{\K}$ has an involution of genus $g$ then $\pi\geq 2g-1$ and the quadrics of
$\P(H^0(I_{C_{\K}}(2)))$ are a base-harmonic system respect to the base spaces $\P^{g-1}$,
$\P^{\pi-g-1}$ that contains $(g-1)(\pi-g-2)$ independent base quadrics. Conversely, these
conditions are sufficient to grant the existence of an involution, except when:

\begin{enumerate}

\item $\pi=6,g=2$ and $C_{\K}$ has a $g^2_5$; or

\item $\pi=2g$,$2g+1$ or $2g+2$ and $C_{\K}$ is trigonal.

\end{enumerate}

\end{teo2} \qed

\begin{rem2}
{\em Let us study what happens at the two exceptions:

\begin{enumerate}

\item Suppose that $C_{\K}$ is a canonical curve of genus $6$ with a $g^2_5$, that is,it is
isomorphic to a smooth plane curve of degree $5$. Suppose that the quadrics of
$\P(H^0(I_{C_{\K}}(2)))$ are a base-harmonic system respect to the base spaces $\P^1$ and
$\P^3$. It induces a harmonic involution on the Veronese surface and then, an involution on
the plane.
 Obviously, the generic plane curve of degree $5$ of the plane is not invariant by this
involution. So in this case the hypothesis of the above theorem are not sufficient.

However, there are smooth quintic plane curves invariant by an involution. For example, we can
take the quintic curve $f(x_0,x_1)-x_2^4x_0=0$ on the plane, where $f(x_0,x_1)$ is a generic
homogeneous polynomial of degree $5$. This curve is smooth an it's invariant by the involution
$$
\begin{array}{ccc}
{x_0\lrw x_0;}&{x_1\lrw x_1;}&{x_2\lrw -x_2}\\
\end{array}
$$

\item Now, suppose that $C_{\K}$ is trigonal. Then it lies on a rational
ruled surface $S_e=\P(\Te_{P^1}\oplus \Te_{P^1}(-e))$ in the linear systems $|3X_0+af|$. The
canonical embedding is obtained by the linear system $X_0+(a-e-2)f$ on the ruled surface.

If $\P(H^0(I_{C_{\K}}(2)))$ is a base-harmonic system then it defines a harmonic involution
on the ruled surface $S_e$. Moreover, we have an induced involution in the linear system
$|3X_0+a f|$. The generic curve of this linear system is not invariant by the involution. We
see that the hypothesis of the theorem are not sufficient.

On the other hand, there are smooth curves on these linear systems invariant by the involution.
Let us see an example. Consider the rational ruled surface $S_0\cong \P^1\times \P^1$ with
coordinates $[(x_0,x_1),(y_0,y_1)]$. We can take the curve on $S_0$ with equation:
$$
x_0^ny_0^3-x_0^ny_0y_1^2+x_1^ny_1^3+x_1^ny_0^2y_1=0
$$
with $n\geq 5$ even.
This is a smooth curve of type $(3,n)$ on the linear system $3X_0+nf$. Moreover it is
invariant by the involution
$$
\begin{array}{cccc}
{x_0\lrw x_0;}&{x_1\lrw -x_1;}\\
{y_0\lrw y_0;}&{y_1\lrw -y_1}\\
\end{array}
$$

\end{enumerate} }\qed

\end{rem2}

\end{document}